\documentclass[12pt]{article}

\usepackage{graphics,amsthm, amsmath, amssymb}

\usepackage{psfig}
\usepackage{color}

\expandafter\chardef\csname pre amssym.def at\endcsname=\the\catcode`\@
\catcode`\@=11

\def\undefine#1{\let#1\undefined}
\def\newsymbol#1#2#3#4#5{\let\next@\relax
 \ifnum#2=\@ne\let\next@\msafam@\else
 \ifnum#2=\tw@\let\next@\msbfam@\fi\fi
 \mathchardef#1="#3\next@#4#5}
\def\mathhexbox@#1#2#3{\relax
 \ifmmode\mathpalette{}{\m@th\mathchar"#1#2#3}%
 \else\leavevmode\hbox{$\m@th\mathchar"#1#2#3$}\fi}
\def\hexnumber@#1{\ifcase#1 0\or 1\or 2\or 3\or 4\or 5\or 6\or 7\or 8\or
 9\or A\or B\or C\or D\or E\or F\fi}

\font\tenmsa=msam10
\font\sevenmsa=msam7
\font\fivemsa=msam5
\newfam\msafam
\textfont\msafam=\tenmsa
\scriptfont\msafam=\sevenmsa
\scriptscriptfont\msafam=\fivemsa
\edef\msafam@{\hexnumber@\msafam}
\mathchardef\dabar@"0\msafam@39
\def\dashrightarrow{\mathrel{\dabar@\dabar@\mathchar"0\msafam@4B}}
\def\dashleftarrow{\mathrel{\mathchar"0\msafam@4C\dabar@\dabar@}}

\font\tenmsb=msbm10
\font\sevenmsb=msbm7
\font\fivemsb=msbm5
\newfam\msbfam
\textfont\msbfam=\tenmsb
\scriptfont\msbfam=\sevenmsb
\scriptscriptfont\msbfam=\fivemsb
\edef\msbfam@{\hexnumber@\msbfam}
\def\Bbb#1{\fam\msbfam\relax#1}

\theoremstyle{plain}
\newtheorem{theorem}{Theorem}[section]
\newtheorem{lemma}[theorem]{Lemma}

\newtheorem{corollary}[theorem]{Corollary}

\theoremstyle{definition}

\newtheorem{example}[theorem]{Example}

\newtheorem{remark}[theorem]{Remark}

\def\para{\vspace{7mm}}

\def\res{{\rm Res}}
\def\deg{{\rm deg}}

\def\cP{{\mathcal P}}
\def\cQ{{\mathcal Q}}
\def\dd{{\mathcal D}}
\def\cc{{\mathcal C}}
\def\ccc{\overline{\mathcal C}}
\def\ee{{\mathcal E}}

\def\N{{\mathbb N}}

\def\qed{\hfill  \framebox(5,5){}}

\def\para{\vspace{2 mm}}

\def\card{{\rm Card}}

\def\R{{\mathbb R}}
\def\C{{\mathbb C}}

\def\K{{\mathbb K}}
\def\F{{\mathbb F}}

\def\L{{\mathbb L}}
\def\Q{{\mathbb Q}}

\def\proj{{\mathbb P}}

\def\card{{\rm card}}
\def\tdeg{{\rm tdeg}}

\def\sup{{\rm sup}}
\def\inf{{\rm inf}}

\def\lim{{\rm lim}}

\def\tgh{ \mathbb{T}(P,\ee^h)}

\def\ddd{{\rm d}}
\def\HH{{\rm H}}

\def\deg{{\rm deg}}

\def\cP{{\mathcal P}}
\def\cE{{\mathcal E}}

\def\cT{{\mathcal T}}

\def\cQ{{\mathcal Q}}
\def\cc{{\mathcal C}}

\begin{document}

\title{Approximate Parametrization of Space Algebraic Curves\thanks{
This work has been developed, and partially supported, by the
Spanish "Ministerio de Ciencia e Innovaci\'on" under the Project
MTM2008-04699-C03-01, and by the "Ministerio de Econom\'{\i}a y
Competitividad" under the project MTM2011-25816-C02-01.  All authors
belong to  the Research Group ASYNACS (Ref. CCEE2011/R34).}}


\author{
Sonia L. Rueda  \\
Dpto. de Matem\'atica Aplicada \\
      E.T.S. Arquitectura,
 Universidad Polit\'ecnica de Madrid \\
        E-28040 Madrid, Spain \\sonialuisa.rueda@upm.es
\and Juana  Sendra \\
Dpto. de Matem\'atica Aplicada a la I.T. de Telecomunicaci\'on  \\
       E.U.I.T.Telecomunicaci\'on, Universidad Polit\'ecnica de Madrid \\
        E-28031 Madrid, Spain \\jsendra@euitt.upm.es
\and
 J. Rafael  Sendra\\
Dpto. de Matem\'aticas \\
        Universidad de Alcal\'a \\
      E-28871 Madrid, Spain  \\
rafael.sendra@uah.es
}
\date{}          
\maketitle

\begin{abstract}
Given a non-rational real space curve and a tolerance $\epsilon>0$, we present an algorithm to approximately parametrize the curve. The algorithm checks whether a planar projection of the space curve is $\epsilon$-rational and, in the affirmative case, generates a planar parametrization that is lifted to an space parametrization. This output rational space curve is of the same degree as the input curve, both have the same structure at infinity, and the Hausdorff distance between them is always finite.
\end{abstract}

\section*{Introduction}

The development of approximate algorithms for algebraic and geometric problems is an active research area (see e.g. \cite{Robiano} and \cite{Stetter-book}) that focuses on different problems as, for instance,  the computation of gcds (see \cite{cor:gi:trwat}, \cite{em:ga:lo}, \cite{Pan1}), the factorization of polynomials (\cite{Cor2}, \cite{GaRu}, \cite{Kaltofen}), the implicicitization of surfaces (\cite{Cor4}, \cite{Dokken}),  the parametrization of curves and surfaces (\cite{BaRo}, \cite{Josef02}, \cite{Hart}, \cite{PSS}, \cite{PSS2}, \cite{PRSS}), etc. These approximate algorithms are applicable by themselves since they face symbolic computation to real world problems. Moreover, those of geometric nature are of special interest in the field of CAGD. For instance, providing parametric representations of algebraic geometric objects helps in some CAGD constructions as surface–-surface intersection or computation of planar sections  (see e.g. Example 5.3. in \cite{PRSS}).

\para

In this context, when  using the term "approximate", there is certain risk of ambiguity since it can have a double meaning (see e.g. \cite{Stetter}). Let us clarify what we mean in our case. Here,  as in our previous papers \cite{PSS}, \cite{PSS2}, \cite{PRSS}, with the term "approximate"
 we do not mean "numerical" but something different (see introduction to \cite{PRSS} for further details): our input is the perturbation of an unknown input; once the perturbation is received we treat it exactly to  provide  an output  that is close  (in this sense it is "approximate"), under certain distance, to the  theoretical output of the unperturbed  unknown input.

 \para

 More precisely, the problem in this paper is as  follows. Let $\cc^*$ be a rational real space curve defined as the complex-zero set of a finite set ${\cal M }\subset \R[x,y,z]$ of real polynomials; in practice $\card({\cal M})=2$. Nevertheless, instead of getting ${\cal M}$ as input of our problem, we get a new finite subset  ${\cal F}\subset \R[x,y,z]$ (which is a perturbation of $\cal M$) of real polynomials
 that defines a new curve $\cc$,  obviously different to $\cc^*$. Since the genus of a curve is unstable under perturbations, the input curve $\cc$ will have positive genus and hence it will  not be parametrizable by rational functions.
Ideally, the problem would consists in finding the initial curve $\cc^*$ or, even better, a rational parametrization of it. However, this goal is unrealistic. Instead, one might require to find a rational parametrization of the closest rational curve to $\cc$  under certain distance; say, under the Hausdorff distance. Nevertheless, in this paper we deal with a weaker  statement of the problem. Namely, finding a rational parametrization of one rational space curve being close (in comparison to a given tolerance $\epsilon>0$)  to $\cc$ under the Hausdorff distance. Our statement, although may not yield to the best output rational curve,  generates one good answer. This can be seen as a first step for the harder, more general, and theoretical problem of finding the best (in the sense of closest) rational curve being, our solution, meanwhile ready to be used in applications.

\para

 In  \cite{PSS}, the authors show how to solve the problem for the particular case of  $\epsilon$-monomial plane curves (i.e. plane curves having an $\epsilon$-singularity of maximum $\epsilon$-multiplicity); see also \cite{PSS2} for the case of surfaces. Later, in \cite{PRSS}, the problem was solved for the more general case of $\epsilon$-rational plane curves. The current paper is, therefore, the natural continuation of this research
 since it deals with the next step, namely the case of  space curves.

\para

In the unperturbed case, the problem can be solved by birationally projecting the space curve on a plane, checking the genus of the projected curve and, in case of genus zero, parametrizing the plane curve to afterwards inverting the parametrization to a rational parametrization of the input curve (see e.g. \cite{SWP} for further details). Now, the situation is more complicated. More precisely the strategy (see box below) is as follows.
$$
{\small \fboxrule 2pt{\fbox{$\begin{array}{ccccc}
\fboxrule 1pt{\fbox{\parbox{2cm}{Let   $\cc\subset \R^3$ } }} & & & & \fboxrule 1pt{\fbox{\parbox{4.5cm}{We get $\overline{\cc}$, parametrized by $\cP(t)=(\frac{p_1(t)}{q(t)},\frac{p_2(t)}{q(t)},\frac{p_3(t)}{q(t)})$, and close to $\cc$.
} }}
\\ $\vector(0,-1){20}$ \vspace*{2mm} && \fboxrule 2pt{\fbox{\parbox{2.8cm}{{\sf General Strategy}} }} & &  \begin{array}{c} \\ \vector(0,1){20}\end{array}
\\
\fboxrule 1pt{\fbox{\parbox{3.5cm}{under a suitable projection (say $\pi_z$ onto the plane $z=0$) we project $\cc$.} }} & & \mathbf{\circlearrowleft} & &\fboxrule 1pt{\fbox{\parbox{4cm}{Applying Chinese Remainder techniques, we lift $\cal D$.} }}
\\
 $\vector(0,-1){20}$  & & & &$\vector(0,1){20}$
\\
\fboxrule 1pt{\fbox{\parbox{3.5cm}{We get a plane curve $\pi_z(\cc)$ satisfying the hypotheses in  \cite{PRSS}.} }} & $\vector(1,0){12}$  & \fboxrule 1pt{\fbox{\parbox{3cm}{We compute the $\epsilon$-genus of $\pi_z(\cc)$; say it is zero. } }} & $\vector(1,0){12}$  & \fboxrule 1pt{\fbox{\parbox{4.5cm}{Algorithm in \cite{PRSS} outputs a rational curve $\cal D$, parametrized as $\cQ(t)=(\frac{p_1(t)}{q(t)},\frac{p_2(t)}{q(t)})$, and   close to $\pi_z(\cc)$.} }}
\end{array}$}}}
$$

We assume some conditions on the original space curve $\cc$ (see  Section \ref{sec-general-assumptions})  such that when it is projected onto a plane  we get a curve satisfying the hypotheses required by the algorithm in \cite{PRSS}; let us denote by $\pi_z(\cc)$ the projected curve, so we are assuming w.l.o.g. in this explanation that projection has been performed on the plane $z=0$. Then, the algorithm in \cite{PRSS} determines whether $\pi_z(\cc)$ is $\epsilon$-rational, where $\epsilon$ is a fixed given tolerance. If $\pi_z(\cc)$ is not $\epsilon$-rational one may try a different projection, but here we simply ask the algorithm to terminate since, although in some examples this seems to work, we do not have any theoretical argumentation to ensure when that projection exists. Otherwise, the algorithm in \cite{PRSS} goes ahead and computes a rational parametrization $\cQ(t)$ of a  plane rational curve $\cal D$. The last step consists in lifting $\cal D$  to a rational space curve $\overline{\cc}$ being close to $\cc$. For this purpose, we first realize that a sufficient condition for the finite Hausdorff distance  requirement, between both curves, is given by the structure at infinity of the input curve $\cc$. Taking into account this fact, and using a Chinese-remainder type  interpolation, we get a rational parametrization of $\overline{\cc}$. As a consequence of this process, we get a rational space curve $\overline{\cc}$ of the same degree as $\cc$, having the same structure at infinity as $\cc$, and such that the Hausdorff distance between $\overline{\cc}$ and $\cc$ is finite.

\para

The structure of the paper is as follows. In Section \ref{sec-general-assumptions}, we introduce the notation that will be used throughout the paper as well as the general assumptions.  Moreover, we comment on the reasons for the inclusion of these assumptions, we discuss how to check them algorithmically, and we show that they (the assumptions) are general enough. Section \ref{sec-projected-curve} is devoted to the projected curve $\pi_z(\cc)$ and, more precisely, to prove that under the general assumptions imposed in \ref{sec-general-assumptions} $\pi_z(\cc)$ satisfies all requirements in the algorithm in \cite{PRSS}. Section \ref{sec-lifted-curve} focuses on how to lift the rational plane curve $\cal D$ (generated by applying the algorithm in \cite{PRSS} to $\pi_z(\cc)$) to the curve $\overline{\cc}$ such that both curves, $\cc$ and $\overline{\cc}$, have the same structure at infinity; note that $\pi_z$ is a birational map between $\cc$ and $\pi_z(\cc)$, but we are lifting ${\cal D}\neq \pi_z(\cc)$. In Section \ref{sec-algorithm-and-examples} we summarize these ideas to derive an algorithm that is illustrated by two examples.
In  Section \ref{sec-error-analysis}, we prove that the Hausdorff distance between the input and output curves, of our algorithm, is always finite.
For this purpose, we briefly study the asymptotes of space curves.

\section{General Assumptions and Notation}\label{sec-general-assumptions}
We consider a computable subfield $\K$ of $\R$, as well as its algebraic closure $\F$; in practice, we may think that $\K=\Q$.
We denote by $\F^2$ and $\F^3$ the affine plane an affine space over $\F$, respectively. Similarly, we denote by $\proj^2(\F)$ and
$\proj^3(\F)$ the projective plane and projective space over $\F$, respectively. Furthermore, if ${\cal A}\subset \F^3$
(similarly if ${\cal A}\subset \F^2$) we denote by ${\cal A}^{*}\subset \F^3$ its Zariski closure, and by
${\cal A}^h\subset \proj^3(\F)$ the projective closure of ${\cal A}^*$. We will consider $(x,y,z)$ as  affine
coordinates and $(x:y:z:w)$ as  projective coordinates. Also, for ${\cal A}$ as above, we denote by
${\cal A}^{\infty}$ the intersection of ${\cal A}^{h}$ with the projective plane of equation $w=0$. In addition, for every polynomial  $H\in \K[x,y,z]$ we denote by $H^{h}(x,y,z,w)$ the homogenization of $H$.

\para

Our method will be based on the projection of  the space curve on a plane.
 Without loss of generality (see below) we will consider
that $z=0$ is the projection plane.   So we introduce the map  $$\pi_z: \F^3 \rightarrow \F^2, (x,y,z)\mapsto (x,y), $$
as well as
$$ \pi_{z}^{h}: \proj^3(\F)\setminus\{ (0:0:1:0)\} \rightarrow \proj^2(\F), (x:y:z:w)\mapsto (x:y:w).
$$
 Our main object of study will  be an {\sf irreducible} (over $\F$)  affine {\sf real} ({\sf non-planar}) space curve $\cc \subset \F^3$. Although, in practice, in most cases,  $\cc$ will be given
  by two generators, we   present the results for the general case  where  a finite set of generators is provided.
  Therefore, we assume that  $\cc$ is given as the zero-set (over $\F$) of a finite set of real polynomials
  ${\cal F}=\{F_1,\ldots,F_s\}\subset\K[x,y,z]$, $s\geq 2$.  $\epsilon$ will be the tolerance and we assume that $0< \epsilon <1$. In addition, we assume the following:

\para

\noindent \underline{{\sf General Assumptions}}
\begin{enumerate}
\item The cardinality of $\cc^\infty$ is $\deg(\cc)$.
\item $\pi_z: \cc \rightarrow \pi_{z}(\cc)^*$ is birational and $\deg(\cc)=\deg(\pi_z(\cc)^*)$.
\item $(1:0:\lambda:0), (0:1:\mu:0), (0 : 0 : 1 : 0)\not\in \cc^\infty$ for any $\lambda,\mu \in \F$.
\item If $(1:\lambda:\mu:0), (1:\lambda:\mu^*:0)\in \cc^\infty$ then $\mu= \mu^*$.
\item The coefficient of $F_1$ in $z^{\tdeg(F_1)}$ is a non-zero constant; where $\tdeg$ denotes the total degree of $F_1$.
\end{enumerate}

We briefly comment on the reasons for the inclusion of the above assumptions,  and we describe how to check them algorithmically.
The condition on irreducibility is natural since rational curves are irreducible varieties. In any case,   one can always consider
the irreducible decomposition of the input to apply the results to each of the irreducible components. The assumption on
the reality of the curve is included because of the nature of the problem, but the theory can be similarly developed for
the case of complex non-real curves. The exclusion of  planar curves is to  simplify the exposition.
Note that this is not a loss of generality, since one can always apply the algorithm in  \cite{PRSS}.

Concerning the general assumptions, condition (1) will play a fundamental role in the error analysis, and it will be used to  ensure that the Hausdorff distance between output and input is always finite. The birational requirement in condition (2) is introduced to reduce the problem to a
plane curve after projection, and the degree fact will be used, in combination with (3) and (4), to ensure that $\pi_z(\cc)^*$ has as many different points at infinity
as its degree; this condition is required by the algorithm in \cite{PRSS}.
Conditions (3) and (4) are related to the projection $\pi_z$. On one hand,  $(1:0:\lambda:0), (0:1:\mu:0) \not\in \cc^\infty$, for any $\lambda,\mu \in \F$, ensures that $(1:0:0),(0:1:0)\not\in\pi_z(\cc)^\infty$ which is a requirement for the algorithm in
\cite{PRSS} to be applied to $\pi_z(\cc)^*$. On the other, $(0 : 0 : 1 : 0)\not\in \cc^\infty$ guarantees that $\pi_{z}^{h}$ is well defined on $\cc^h$. In addition, conditions (3)-(4) ensure that $\pi_{z}^{h}$ is injective on $\cc^\infty$. Condition (5) is also introduced to guarantee that $\pi_z(\cc)^*$ satisfies the hypotheses in \cite{PRSS}  (see Theorem \ref{theorem-curva-proyectada}). Note that this property can always be achieved by means of  a suitable orthogonal affine change of coordinates, and hence preserving distances.

\para

Taking into account that, in practice, $\cc$ is expected to
 come from the perturbation of a
rational space real curve, in general, all above conditions will hold.
Nevertheless, let us discuss how to decide algorithmically whether a given input satisfies them.
Checking the irreducibility of $\cc$ can be approached by checking whether the corresponding ideal is prime
(see, for instance, Section 4.5 in \cite{Cox:ideals} or \cite{Gianni-Trager-Zacharias}).
In order to check the reality, one can apply cylindrical algebraic decomposition techniques to decide the
existence of real regular points (see e.g. \cite{Johnson}).  The non-planarity of $\cc$ can be deduced from a
Gr\"obner basis of $\cal F$. Let us now deal with the general assumptions. One
can compute $\cc^h$  by homogenizing a Gr\"obner basis of $\cal F$, w.r.t. a graded order (see e.g. page  382 in \cite{Cox:ideals}). The
degree of $\cc$ can de determined by counting the number of intersections of $\cc$ with a generic plane; in fact a randomly chosen plane
might be enough. So, (1) can also be checked. Also, (3) and (4) are checkable. Condition (2) can be analyzed by direct application of elimination theory techniques. Condition (5) is trivially checkable.

\para

In the previous description, we have considered that $z=0$ is the projection plane. Indeed,
conditions (3)-(5)  depend on this fact. So, if any of these conditions fails we might either
consider a suitable orthogonal affine change of coordinates or choose another projection plane.
Also, if (2) fails we need to find a different projection plane. We recall that, for  almost
every plane, the corresponding projection is birational over $\cc$ and that for almost every
plane the number of  intersection points of the plane with  $\cc$ is $\deg(\cc)$. Therefore, the combination of these two facts with
 Lemma \ref{projection-degree} ensures that condition (2) must be achieved by taking the projection plane randomly.

\para

\begin{lemma}\label{projection-degree}
Let  $\Pi\subset \F^3$ be a plane  such that  $\card(\cc\cap \Pi)=\deg(\cc)$,  let $u$ be a (non-zero) parallel   vector to $\Pi$ and non-parallel
to the vectors in $\{ P-Q\,|\, P,Q\in   \cc\cap \Pi, \,P\neq Q\}$, and let $\Pi^{u}$ be any plane orthogonal to $u$.
 Then, $\deg(\pi_{\Pi^{u}}(\cc)^*)=\deg(\cc),$ where $\pi_{\Pi^{u}}$ is the projection map from $\F^3$ onto $\Pi^u$.
\end{lemma}

\para

\noindent {\bf Proof.} Let $d=\deg(\cc)$ and  $\cc\cap \Pi=\{P_1,\ldots,P_d\}$, and let $L$ be the line $\Pi\cap \Pi^u$. By construction, $\{\pi_{\Pi^{u}}(P_i)\}_{i=1,\ldots,d}\subset
\pi_{\Pi^{u}}(\cc)^*\cap L$. Since $u$ is not parallel to $P_i-P_j$, with $i\neq j$, then $\card(\{\pi_{\Pi^{u}}(P_i)\}_{i=1,\ldots,d})=d$. Therefore,
$\deg(\cc) \leq \deg(\pi_{\Pi^{u}}(\cc)^*)$. Now, let $L'$ be a line in $\Pi^u$ such that $L'\cap (\pi_{\Pi^{u}}(\cc)^*\setminus \pi_{\Pi^{u}}(\cc))=\emptyset$
and such that $\card(L' \cap \pi_{\Pi^{u}}(\cc)^*)=\deg(\pi_{\Pi^{u}}(\cc)^*)$; note that almost all lines in $\Pi^u$ satisfy this property. Let
$L' \cap \pi_{\Pi^{u}}(\cc)^*=\{Q_1,\ldots,Q_{d'}\}$ and let $\Pi'$ be the plane containing $L'$ and being parallel to $u$; note that $u$ is normal to $\Pi^{u}$, and hence $L'$ is not parallel to $u$. Because of the construction $\pi_{\Pi^{u}}^{-1}(Q_i)\cap \cc \neq\emptyset$ and it is contained in $\Pi'$. Therefore, $\cup_{i=1}^{d'} \pi_{\Pi^{u}}^{-1}(Q_i)\cap \cc\cap \Pi'$ has cardinality at least $d'$, and hence $\deg(\pi_{\Pi^{u}}(\cc)^*)=d'\leq \deg(\cc)$. \qed

\section{The Projected Curve}\label{sec-projected-curve}

In this section, we analyze the basic properties of the projected curve $\pi_z(\cc)^*$. In particular, we show that it satisfies the hypotheses in \cite{PRSS}. We  recall that, since $\pi_z$ is birational and $\cc$ is irreducible, $\pi_z(\cc)^*$ is irreducible.  We start  with a technical lemma on Gr\"obner bases


\para

\begin{lemma}
\label{lemma-GB} Let ${\Bbb L}\subset \F$ be a field and $G_1,\ldots,G_m\in {\Bbb L}[x_1,\ldots,x_n]$ be such that $\deg_{x_{n}}(G_1)=\tdeg(G_1)>0$; where $\tdeg$ denotes the total degree.
Let $\{H_1,\ldots,H_r\}$ be a Gr\"obner basis of $(G_1,\ldots,G_m)$ w.r.t. the graded lex order with $x_1<\cdots<x_n$. Then, there exists $i\in \{1,\ldots,r\}$ such that $\deg_{x_n}(H_i)=\tdeg(H_i)$. Moreover, if the variety defined by $\{G_1,\ldots,G_m\}$ over $\F$ is not empty then $\deg_{x_n}(H_i)=\tdeg(H_i)>0$.
\end{lemma}

\para

\noindent {\bf Proof.}
Let $\ell=\deg_{x_n}(G_1)$. Since $\ell=\tdeg(G_1)$, because of the ordering, the leading term of $G_1$ is $x_{n}^{\ell}$. Now, by Exercise 5, page 78, \cite{Cox:ideals}, there exists $i\in \{1,\ldots,r\}$ such that the leading term of $F_i$ divides $x_{n}^{\ell}$. Finally, because of the ordering, $\tdeg(H_i)=\deg_{x_n}(H_i)$.  Now, if the variety of $\{G_1,\ldots,G_m)\}$ is empty, we assume that the Gr\"obner basis is normal (this does not affect to the previous reasoning). By Theorem 8.4.3 in \cite{Winkler}, $H_i$ is not constant. So $\tdeg(H_i)>0$. \qed

\para
The next lemma shows how generalized resultants can be used to compute the projection.

\para

\begin{lemma}\label{lemma-resultante-generalizada}
Let
$$ F_\Delta(x,y,z,\Delta)=\left\{\begin{array}{lll} F_2+\Delta F_3+\cdots + \Delta^{s-2} F_s & &\mbox{if $s> 2$} \\
F_2 & & \mbox{if $s=2$} \end{array} \right.,$$ where $\Delta$ is a
new variable, and let
 $F_{\Delta}^{h}(x,y,z,w,\Delta)$ be the homogenization of $F_\Delta(x,y,z,\Delta)$ as a
  polynomial in $\K[\Delta][x,y,z]$; recall that $w$ is the variable used for homogenization. Let
$$ R=\res_z(F_1,F_\Delta)= \sum_{j=0}^{m} \alpha_j(x,y) \Delta^j,\,\, S=\res_{z}(F_{1}^{h},F_{\Delta}^{h})=\sum_{i=0}^{m'} \beta_{i}(x,y,w) \Delta^i. $$
It holds that
 \begin{enumerate}
 \item  $\pi_z(\cc)^*$  is the affine plane curve defined by
$  \gcd(\alpha_0,\ldots,\alpha_m)$, and $m=m'$.
\item If ${\cal F}=\{F_1,\ldots,F_s\}$ is a Gr\"obner basis of the ideal generated by $\cal F$, w.r.t. the graded lex order with $x<y<z$, then $\pi_z(\cc)^h$ is the projective plane curve defined by
$\gcd(\beta_0,\ldots,\beta_{m})$.

\end{enumerate}
\end{lemma}
\noindent {\bf Proof.} We start recalling that, in general,
$\pi_z(\cc)$ is not closed and hence one may need to distinguish
through the proof between $\pi_z(\cc)$ and its Zariski closure
$\pi_{z}(\cc)^*$.

\noindent (1) We first prove that $\pi_z(\cc)^*$ is the variety
defined by $\{\alpha_0,\ldots,\alpha_m\}$. Indeed, let $(a,b)\in
\pi_z(\cc)$. Then, there exists $c\in \F$ such that $P=(a,b,c)\in
\cc$.  So, $F_1(P)=0, F_\Delta(P,\Delta)=0$. Thus,
$R(a,b,\Delta)=0$, and hence $
\alpha_0(a,b)=\cdots=\alpha_{m}(a,b)=0$. Conversely, let
$\alpha_0,\ldots,\alpha_m$ vanish at $(a,b)$. Then
$R(a,b,\Delta)=0$. Now, since
     $\deg_z(F_1)=\tdeg(F_1)$, there exists $c$ in the algebraic closure of ${\F(\Delta)}$ (indeed $c\in \F$) such that $F_1(a,b,c)=0,F_\Delta(a,b,c,\Delta)=0$.
     Since  $c\in \F$, from $F_\Delta(a,b,c,\Delta)=0$, we get that  $F_2(a,b,c)=\cdots=F_s(a,b,c)=0$. So, $(a,b,c)\in \cc$ and  $(a,b)\in \pi_z(\cc)$.

     Let $\alpha=\gcd(\alpha_0,\ldots,\alpha_m)$ and $\overline{\alpha}_i$ be such that $\alpha_i=\overline{\alpha}_i \alpha$. Let $\cal V$ and $\cal W$ be the varieties defined by $\alpha$ and $\{\overline{\alpha}_0,\ldots,\overline{\alpha}_m\}$, respectively. Then, $\pi_z(\cc)^*={\cal V}\cup {\cal W}$. Since $\cc$ is irreducible and $\pi_z$ birational, we have that $\pi_z(\cc)^*$ is an irreducible curve. So, $\cal V$ is 1-dimensional and $\cal W$ is either empty or 0-dimensional. In any case, because of the irreducibility, ${\cal W}\subset \cal V$. So $\pi_z(\cc)^*=\cal V$.

  Finally, let us see that $m=m'$. We assume w.l.o.g. that all $\beta_i$ are non-zero. Since   $\deg_z(F_1)=\deg_z(F_{1}^{h})$,  taking into account how the resultants behave under specializations (see, e.g. Lemma 4.3.1. in \cite{Winkler}), we get that $R(x,y,\Delta)=S(x,y,1,\Delta)$, up to multiplication by a non-zero constant. Moreover,  $S$ is homogeneous as a polynomial in $\F[\Delta][x,y,w]$. Thus, $\beta_i$ are homogeneous (of the same degree). Therefore, $\beta_i(x,y,1)$ does not vanish. So, $m=\deg_\Delta(R)=\deg_\Delta(S(x,y,1,\Delta))=\deg_\Delta(S)=m'$.

\para

\noindent (2) We first prove that $w$ does not divide $S$. Let $S=w
M(x,y,w,\Delta)$. Then, for all $(a,b)\in \F^2$, since
$\deg(F_{1}^{h})=\deg_z(F_{1}^{h})$, there exists $c$  in the
algebraic closure of  ${\F(\Delta)}$
 such that $F_{1}^{h}(a,b,c,0)=F_{\Delta}^{h}(a,b,c,0,\Delta)=0$. Therefore, since indeed $c\in \F$,  $F_{i}^{h}(a,b,c,0)=0, i=1,\ldots,s$;  let us call $\rho(a,b)$ the corresponding $c$ associated to $(a,b)$.
Then, the infinitely many points   $\{ (1: n :\rho(1,n):0) \}_{n\in
{\mathbb{N}}}$ are included in the intersection of $\cc^h$  with the
plane $w=0$, which is a contradiction.

From $R(x,y,\Delta)=S(x,y,1,\Delta)$ we get that
$\alpha_j(x,y)=\beta_j(x,y,1)$. Therefore,  $\alpha_{j}^{h}
w^{n_j}=\beta_j$, for some $n_j\in {\mathbb{N}}$. Moreover, since
$w$ does not divide $S$, there exists $i_0\in \{0,\ldots,m\}$ such
that $\alpha_{i_0}^{h}=\beta_{i_0}$ and
$\gcd(\alpha_{i_0}^{h},w)=1$.

Let $\alpha=\gcd(\alpha_0,\ldots,\alpha_m)$ and
$\gamma=\gcd(\alpha_{0}^{h},\ldots,\alpha_{m}^{h})$. We see that
$\alpha^{h}=\gamma$. Let $\alpha_{i}=\alpha \overline{\alpha}_i$.
Then, $\alpha_{i}^{h}=\alpha^{h} \overline{\alpha}_{i}^{h}$. So,
$\alpha^h$ divides $\gamma$. Conversely, let $\alpha_{i}^{h}=\gamma
\tilde{\alpha}_i$. Then,  $\gamma(x,y,1)$ divides
$\alpha_{i}^{h}(x,y,1)=\alpha_i$. Therefore, $\gamma(x,y,1)$ divides
$\alpha$. In addition, since $\alpha^h$ divides $\gamma$, $\alpha$
divides $\gamma(x,y,1)$. Hence, up to multiplication by non-zero
constants, $\alpha=\gamma(x,y,1)$. Therefore, since  $w$ does not
divide $\gamma$, we get that $\alpha^h=\gamma$.

Finally, it remains to prove that
$\gamma=\gcd(\beta_0,\ldots,\beta_m)$. We know that
$\gcd(\beta_0,\ldots,\beta_m)=\gcd(\alpha_{0}^{h}
w^{n_0},\ldots,\alpha_{m}^{h}w^{n_m})$. Let $a=\gcd(\alpha_{0}^{h}
w^{n_0},\ldots,\alpha_{m}^{h}w^{n_m})$. Clearly $\gamma$ divides
$a$. Conversely, $a$ divides $\alpha_{i_0}^{h}$ (see above). Since
$\gcd(\alpha_{i_0}^{h},w)=1$, then $\gcd(a,w)=1$. Therefore, $a$
must divide all $\alpha_{j}^{h}$. Hence, $a$ divides $\gamma $.

Summarizing
$\gcd(\alpha_0,\ldots,\alpha_m)^h=\gcd(\alpha_{0}^{h},\ldots,\alpha_{m}^{h})=\gcd(\alpha_{0}^{h}
w^{n_0},\ldots,\alpha_{m}^{h}w^{n_m})
=\gcd(\beta_0,\ldots,\beta_m)$. \qed

\para

\para
\begin{remark}\label{remark-resultant}
We observe that in the proof of Lemma \ref{lemma-resultante-generalizada}, from all the hypotheses imposed in Section \ref{sec-general-assumptions}, we have  only used the following: general assumption (5) is used in both (1) and (2). The fact that $\cc$ has dimension 1 and that is irreducible is used in (1), jointly with the fact that $\pi_z$ is finite. Finally, in (2), we use that $\cc^h$ intersects the plane $w=0$ in finitely many points.

%
%
\end{remark}

\para

\para

We finish the section by stating the main properties of the projected curve.

\para

\begin{theorem}\label{theorem-curva-proyectada}
It holds that
\begin{enumerate}
\item $\pi_{z}^{h}(\cc^\infty)=\pi_{z}(\cc)^\infty$.
\item $\card(\pi_{z}(\cc)^\infty)=\deg(\pi_{z}(\cc)^*)$.
\item $(1:0:0),(0:1:0)\not\in \pi_{z}(\cc)^\infty$.
\end{enumerate}
\end{theorem}
\para

\noindent {\bf Proof.}
(1) Because of Lemma \ref{lemma-GB}, we can assume w.l.o.g. that $\{F_1,\ldots,F_s\}$ is a Gr\"obner basis w.r.t the
 graded lex order with $x<y<z$, and that $\deg_z(F_1)=\tdeg(F_1)$. Also, let $F_{\Delta}, F_{\Delta}^{h}, S,R, \alpha_i, \beta_i$ be as in Lemma \ref{lemma-resultante-generalizada}, and $\alpha=\gcd(\alpha_0,\ldots,\alpha_m)$.

Note that $\cc^{\infty}$ is the zero set in ${\Bbb P}^3(\F)$ of $\{F_{1}^{h}(x,y,z,0),\ldots,F_{s}^{h}(x,y,z,0)\}$. So, since  $\cc^{\infty}$ is zero-dimensional, then $\gcd(F_{1}^{h}(x,y,z,0),\ldots,F_{s}^{h}(x,y,z,0))=1$. In addition, by Lemma \ref{lemma-resultante-generalizada},
$\pi_{z}(\cc)^{\infty}$ is the zero set in ${\Bbb P}^2(\F)$ of $\alpha^{h}(x,y,0)$.

 Now, let $(a:b:c:0)\in \cc^{\infty}$. Then, $F_{1}^{h}(a,b,c,0)=F_{\Delta}^{h}(a,b,c,0)=0$. Therefore, $S(a,b,0,\Delta)=0$. Thus, $\beta_{i}(a,b,0)=0$. By Lemma \ref{lemma-resultante-generalizada}, we know that $\alpha^h(x,y,w)=\gcd(\beta_0,\ldots,\beta_m)$. So, $\beta_i=\alpha^h\overline{\beta}_i$. Let us assume that $\alpha^h(a,b,0)\neq 0$ (i.e. that $(a:b:0)\not\in \pi_z(\cc)^{\infty}$).  By general assumption (3),   $a,b$ cannot be both zero.
 We   assume w.l.o.g. that $a=1$. Then, $\overline{\beta}_i(1,b,0)=0$ for all $i$. We  now consider the polynomials $H_i(y,z,w)=F_{i}^{h}(1,y,z,w)$ as well as  the affine variety ${\cal D}$ defined by them. Note that, since $(1:b:c:0)\in \cc^h$,  $(b,c,0)\in{\cal D}\neq \emptyset$. Moreover, since $\cc$ is irreducible, $\cal D$ is an irreducible curve. Furthermore, $\tdeg(H_1)=\deg_z(H_1)=\deg_z(F_1)>0$.  Furthermore, since $\cc$ is not planar, $\cal D$ is not a line perpendicular to the plane $z=0$, so $\pi_z$ is finite over $\cal D$. Also, since $\cc$ is not planar, ${\cal D}^h$ intersection $x=0$ has only finitely many points.
 Furthermore, because of Lemma \ref{lemma-GB}, we can assume w.l.o.g. that $\{H_1,\ldots,H_s\}$ is a Gr\"obner basis w.r.t the
 graded lex order with $w<y<z$, and that $\deg_z(H_1)=\tdeg(H_1)$.
 Thus, $\cal D$ satisfies the hypotheses of Lemma \ref{lemma-resultante-generalizada} (see Remark \ref{remark-resultant}).
 Let $H_{\Delta}$ as in  Lemma \ref{lemma-resultante-generalizada}, let $T(y,w,\Delta)=\res_z(H_1,H_\Delta)$, and let $S,\alpha_{i},\alpha^{h},\beta_i$ be as in the proof of Lemma \ref{lemma-resultante-generalizada}.
   Reasoning as in the proof Lemma \ref{lemma-resultante-generalizada}, we get that $\deg_{\Delta}(T)=\deg_{\Delta}(S)$ and that, if $T=\sum_{i=0}^{m} \rho_{i}\Delta^i$ then  $\rho_i(y,w)=\beta_i(1,y,w)$. Moreover, by Lemma \ref{lemma-resultante-generalizada}, we get that $\pi_z({\cal D})^*$ is defined by $\rho=\gcd(\rho_0,\ldots,\rho_m)$; note that $\pi_z({\cal D})^*$ is irreducible, and hence $\rho$ is an irreducible polynomial.

 From $\rho_i(y,w)=\beta_i(1,y,w)=\alpha^{h}(1,y,w)\overline{\beta_i}(1,y,w)$ we get that $\alpha^{h}(1,y,w)$ divides $\rho(y,w)$.
Also, since $\cc$ is not planar, $\pi_z(\cc)^*$ is not a line, and hence $\alpha^{h}(1,y,w)$ is not constant. Thus, since $\rho$ is irreducible, we get that, up to multiplication by non-zero constants, $\rho(y,w)=\alpha^{h}(1,y,w)$. Finally, from $(b,c,0)\in \cal D$, we get that $\rho(b,0)=\alpha^{h}(1,b,0)=0$, which is a contradiction. This proves that $\pi_{z}^{h}(\cc^\infty)\subset \pi_{z}(\cc)^\infty$.

On the other hand, because of general assumptions (3) and (4) one has that $\pi_{z}^{h}$ is injective over $\cc^\infty$, and hence  we get that $\card(\pi_{z}^{h}(\cc^\infty))=\card(\cc^\infty)$. Then, from general assumptions (1) and (2), we get that $\pi_{z}^{h}(\cc^\infty)= \pi_{z}(\cc)^\infty$.

 \para

 \noindent (2)  Because of general assumptions (3) and (4), $\card(\cc^\infty)=\card(\pi_{z}^{h}(\cc^\infty))$ and, by general assumptions (1) and (2), $\card(\cc^\infty)=\deg(\cc)=\deg(\pi_z(\cc)^*)$. Now the proof ends by applying statement (1) in this theorem.

\para

\noindent (3)  It follows from general assumption (3) and statement (1) in this theorem. \qed

\section{The Lifted Curve}\label{sec-lifted-curve}

In Theorem \ref{theorem-curva-proyectada} we have seen that, under the assumptions introduced in Section \ref{sec-general-assumptions}, $\pi_{z}(\cc)^*$ satisfies the hypotheses required by the parametrization algorithm in \cite{PRSS}.
In this situation, we apply algorithm in \cite{PRSS} to $\pi_{z}(\cc)^*$. If $\pi_{z}(\cc)^*$ is not $\epsilon$-rational, then we can not use $\pi_{z}(\cc)^*$ to parametrize $\cc$
approximately by this method. However,  it might be that there exists another projection such that the projected curve is $\epsilon$-rational and hence the method applicable to this other projection. Nevertheless, we have not researched in this direction leaving this as a future research line. So, let us suppose that $\pi_{z}(\cc)^*$ is $\epsilon$-rational, and let
$$\cQ(t)=\left(\frac{p_1(t)}{q(t)},\frac{p_2(t)}{q(t)}\right) $$
be the parametrization output by the algorithm in \cite{PRSS}. Let $\cal D$ be the rational plane curve parametrized by $\cQ(t)$.
We want to lift $\dd$ from $\F^2$ to a rational curve $\ccc$ in  $\F^3$. For this purpose, in order to guarantee that the Hausdorff distance between $\cc$ and $\ccc$ is finite (see Corollary \ref{cor-distance}), we will associate to $\dd$ a rational curve $\ccc$ in $\F^3$  such that $\pi_{z}(\ccc)=\dd$, $\deg(\ccc)=\deg(\cc)$ and  $\cc^\infty=\ccc^\infty$.

\para

We know that $\pi_z(\cc)^*$ and $\dd$ have the same degree and the same structure at infinity (see Theorem 4.5. in \cite{PRSS}). Thus, by Theorem \ref{theorem-curva-proyectada}, $\dd^\infty=\pi_{z}^{h}(\cc^\infty)$.
 In addition, it also holds that $\deg(p_i)\leq \deg(q)=\deg(\dd)$ (see proof of Lemma 4.2 in \cite{PRSS}). Moreover, by construction, $\gcd(p_i,q)=1$ (see Step 10 in the algorithm in \cite{PRSS}). Furthermore, $(p_1(t):p_2(t):q(t))$ reaches all points in $\dd^\infty$ (see proof of Theorem 4.5. in \cite{PRSS}). Therefore, since $\card(\dd^\infty)=\deg(\dd)$ (see Theorem \ref{theorem-curva-proyectada} and note that $\dd^\infty=\pi_{z}^{h}(\cc^\infty)$), $q(t)$ is square-free. Thus, if $\{\xi_1,\ldots,\xi_d\}$ are the roots of $q(t)$,
 $$\dd^\infty=\left\{\left(1: \frac{p_2(\xi_i)}{p_1(\xi_i)}:0\right)\right\}_{i=1,\ldots,d},$$
 because of general assumption (4), for every $i$ there exists a unique $\chi_i\in \F$ such that
  $$\cc^\infty=\left\{\left(1: \frac{p_2(\xi_i)}{p_1(\xi_i)}:\chi_i:0\right)\right\}_{i=1,\ldots,d}.$$
Note that, if  $\{G_1,\ldots,G_m\}$ is a Gr¨\"obner basis of $\{F_1,\ldots,F_s\}$ w.r.t. the graded lex order with $x<y<z$, then
$\chi_i$ is the root of
\[ \gcd\left(G_{1}^{h}\left(1,\frac{p_2(\xi_i)}{p_1(\xi_i)},z,0\right),\ldots,G_{m}^{h}\left(1,\frac{p_2(\xi_i)}{p_1(\xi_i)},z,0\right)\right). \]
Let $p_3(t)$ be the interpolating polynomial such that $p_3(\xi_i)=p_1(\xi_i)\chi_i$, for $i=1,\ldots,d$; recall that $q(t)$ is square-free. We then define $\ccc$ as the rational curve
\[ \cP(t)=\left(\frac{p_1(t)}{q(t)},\frac{p_2(t)}{q(t)},\frac{p_3(t)}{q(t)}\right). \]
 Note that $\gcd(p_1,p_2,p_3,q)=1$, $q$ is square-free, $\deg(p_1), \deg(p_2)\leq \deg(q)$ and $\deg(p_3)< \deg(q)$.

\para

Taking into account the previous reasonings, we have the following theorem.

\para

\begin{theorem}\label{teorema-lifted-curve}
The lifted curve $\ccc$, defined as above, satisfies that
\begin{enumerate}
\item $\ccc$ is rational.
\item $\cc^\infty=\ccc^\infty$.
\item $\deg(\cc)=\deg(\ccc)$.
\item $\pi_z(\ccc)^*=\dd$.
\end{enumerate}
\end{theorem}

\para

We finish this section  explaining how to compute $\ccc$ (i.e. the polynomial $p_3(t)$) without having to explicitly compute the roots of $q(t)$. The idea is to adapt the Chinese Remainder interpolation techniques. Let $\{G_1,\ldots,G_m\}$ be as above, and let $q(t)=\prod_{j=1}^{\ell} q_j(t)$ be an irreducible factorization of $q(t)$ over $\K$. Now, for each $q_j$ we consider the field $\L=\K(\mu)$, where $\mu$ is the algebraic element over $\K$ defined by $q_j(t)$, as well as the polynomial ring $\L[z]$. Let
$$D_j(z)=\gcd_{\L[z]}\left(G_{1}^{h}\left(1,\frac{p_2(\mu)}{p_1(\mu)},z,0\right),\ldots,G_{m}^{h}\left(1,\frac{p_2(\mu)}{p_1(\mu)},z,0\right)\right),$$
where the gcd is taken in the Euclidean domain $\L[z]$. Because of the previous reasoning, we know that $D_j(z)$ can be expressed as
\[ D_j(z)=(a_j(\mu) z-b_j(\mu))^u\in \L[z], \]
where $u\in \mathbb{N}$. Let $c_j(\mu)$ be the polynomial expression of $b_{j}(\mu)a_{j}(\mu)^{-1}p_1(\mu)$ as an element in $\L$. On the other hand, for $i\neq j$, $\gcd(q_i,q_j)=1$. So, there exist $u_{i,j},u_{j,i}\in \K[t]$ such that $u_{i,j}q_i+u_{j,i}q_j=1$. We introduce the following polynomial
\[ A(t)=\left\{\begin{array}{ll}
c_1(t) \prod_{i=2}^{\ell} u_{i,1}(t)q_{i}(t)+\cdots+c_\ell(t)\prod_{i=1}^{\ell-1} u_{i,\ell}(t)q_{i}(t) & \mbox{if $\ell>1$} \\
c_1(t) & \mbox{if $\ell=1$.}
\end{array} \right. \]
Then, we have the following result.

\para

\begin{lemma}\label{lemma-interpolation}
$p_3(t)$ is the remainder of the division of $A(t)$  by $q(t)$.
\end{lemma}

\para

\noindent {\bf Proof.} If $\ell=1$ the result is trivial. Let $\ell>1$ and let $R(t), Q(t)$ be the remainder and quotient of the division of $A(t)$  by $q(t)$, respectively. Clearly $\deg(R)<\deg(q)$. Now, let $\xi_i$ be a root of $q(t)$; say w.l.o.g. that $\xi_i$ is a root of $q_1(t)$. Then, by construction, $c_1(\xi_i)=\chi_i p_1(\xi_i)$. Therefore,
\[ R(\xi_i)=A(\xi_i)-q(\xi_i)Q(\xi_i)=A(\xi_i)=c_1(\xi_i)u_{2,1}(\xi_i)q_2(\xi_i)\cdots u_{\ell,1}(\xi_i) q_{\ell}(\xi_i).\]
However, for $k\neq 1$, $u_{k,1}(\xi_i)q_k(\xi_i)=1-u_{1,k}(\xi_i)q_1(\xi_i)=1$. So, $R(\xi_i)=c_1(\xi_i)=\chi_i p_1(\xi_i)=p_3(\xi_i)$. \qed

\section{Algorithm and Examples}\label{sec-algorithm-and-examples}

In this section we collect all the ideas developed in the previous sections to derive the approximate parametrization algorithm, and we illustrate it by several examples. For this purpose, we assume that we are given a tolerance $0<\epsilon<1$ as well as an space curve $\cc$ satisfying all the hypotheses imposed in Section \ref{sec-general-assumptions}. Then the algorithm is as follows

\para

\noindent {\sf \underline{Algorithm}}

\para

\begin{enumerate}
\item Compute the defining polynomial  of $\pi_z(\cc)^*$ (apply e.g. Lemma \ref{lemma-resultante-generalizada}).
\item Apply to $\pi_z(\cc)^*$ the parametrization algorithm in \cite{PRSS}. If the plane curve is not $\epsilon$-rational exit returning no parametrization else let $(\frac{p_1(t)}{q(t)},\frac{p_2(t)}{q(t)})$ be the output parametrization.
\item Apply Lemma \ref{lemma-interpolation} to determine $p_3(t)$.
\item Return $(\frac{p_1(t)}{q(t)},\frac{p_2(t)}{q(t)},\frac{p_3(t)}{q(t)})$.
\end{enumerate}

\para

\begin{remark}
Note that, because of Theorem \ref{teorema-lifted-curve}, it holds that the rational curve output by the algorithm has the same degree and structure at infinity as the input curve. Also, as already mentioned in Section \ref{sec-lifted-curve}, if in step 2 we do not  get  $\epsilon$-rationality,
 it does not imply that under another projection one could not  get an $\epsilon$-rational curve. However we can not guarantee  theoretically when such a projection exists.
\end{remark}

\para

In the following examples,  the polynomial $f$ defining $\pi_z(\cc)^*$ and the parametrizations $\cQ(t)$ of $\cal D$,  and $\cP(t)$ of ${\ccc}$,  are expressed with 10-digits floating point coefficients, but the executions have
been performed with  exact arithmetic; the precise data can be found in  \textcolor{blue}{http://www2.uah.es/rsendra/datos.html}.

\para

\begin{example}\label{ex-1}
Let $\cc$ be the space curve defined by the polynomials
\begin{align*}
\mathit{F_1} =&  - {\displaystyle \frac {718945312497}{100}} \,x
 + {\displaystyle \frac {698623125001}{100}} \,y - 671015625\,z
 + 13865578693\,z\,y \\
&\mbox{} - 12118499950\,z\,x + 24392628607\,x\,y - 18401807886\,y
^{2} - 1311877532\,z^{2}\\
\mathit{F_2} =&  - {\displaystyle \frac {431020499999}{25}} \,x +
{\displaystyle \frac {1675347948801}{100}} \,y - 1609143200\,z +
4365980240\,z\,y \\
&\mbox{} - 401217042\,z\,x - 24936051360\,y^{2} - 683547137\,z^{2}
 + 24392628607\,x^{2}
\end{align*}
and let $\epsilon=\frac{1}{100}$. One  can check that $\cc$ satisfies all the hypotheses imposed in Section \ref{sec-general-assumptions}. Moreover, $\deg(\cc)=4$.
The projected curve $\pi_z (\cc)^*$ is defined by the polynomial (see Fig. \ref{fig-1})

\para

\begin{flushleft}
$f(x,y)= { 5.192147942\cdot 10^{29}}\,xy-{ 2.214420657\cdot 10^{28}}\,y-{
 5.059350678\cdot 10^{28}}\,x-{ 2.636990684\cdot 10^{29}}\,{x}^{2}-{
 3.506554787\cdot 10^{42}}\,{x}^{2}y+{ 2.001041491\cdot 10^{42}}\,{y
}^{4}-{ 1.375243688\cdot 10^{42}}\,{y}^{3}+{ 1.181135404\cdot 10^{42
}}\,{x}^{3}+{ 3.822854018\cdot 10^{42}}\,x{y}^{2}-{ 2.315025392
\cdot 10^{40}}\,{y}^{2}-{ 2.990857566\cdot 10^{42}}\,x{y}^{3}-{
 1.221346211\cdot 10^{42}}\,{x}^{2}{y}^{2}+{ 3.915698981\cdot 10^{42
}}\,{x}^{3}y-{ 1.812915331\cdot 10^{42}}\,{x}^{4}.
$
\end{flushleft}



Note that $\deg(\pi_z (\cc)^*)=4$. Applying the
 approximate parametrization algorithm for plane curves in \cite{PRSS} we get that  $\pi_z (\cc)^*$ is $\epsilon$-rational.
  Furthermore  the algorithm outputs the  parametrization $\cQ (t)=(\frac{p_1(t)}{q(t)},\frac{p_2(t)}{q(t)})$
\[ \left({\frac {- 0.4173571408+ 1.171283433\,t- 0.8477221239\,{t}^{2}-
 0.1445883061\,{t}^{3}+ 0.2133409452\,{t}^{4}}{- 0.9059858774+
 1.956830479\,t- 0.6103552658\,{t}^{2}- 1.494650450\,{t}^{3}+{t}^{4}}}
,\right.\]\[\left. {\frac { 0.1828752070\,t+ 0.6268800173\,{t}^{2}- 1.028340444\,{t}^{3}
+ 0.3822448988\,{t}^{4}- 0.1884116000}{- 0.9059858774+ 1.956830479\,t-
 0.6103552658\,{t}^{2}- 1.494650450\,{t}^{3}+{t}^{4}}}  \right) \]

\begin{center}
\begin{figure}[ht]
\centerline{$\begin{array}{cc}
\mbox{\psfig{figure=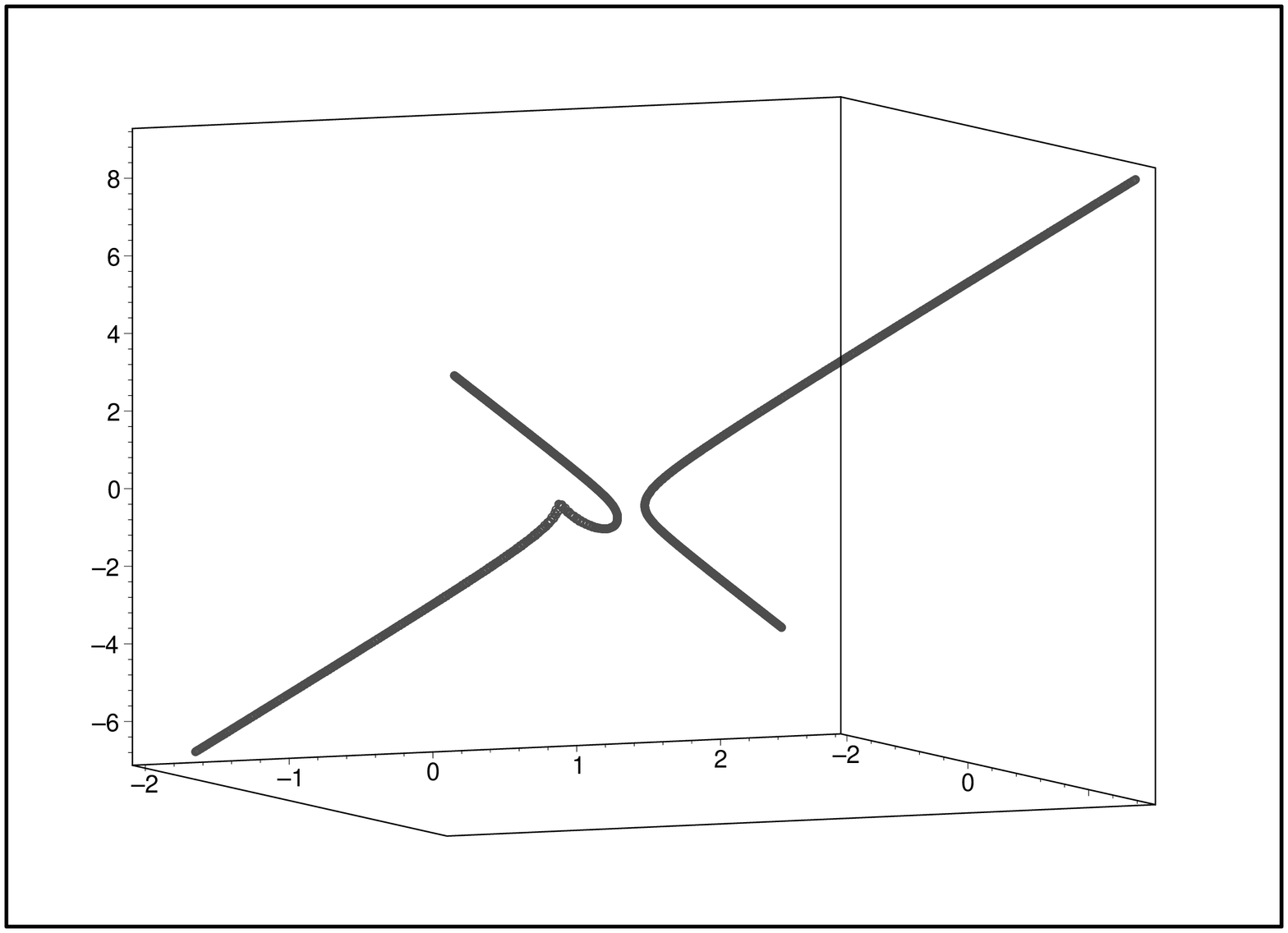,width=6cm,height=6cm,angle=270}} &
\mbox{\psfig{figure=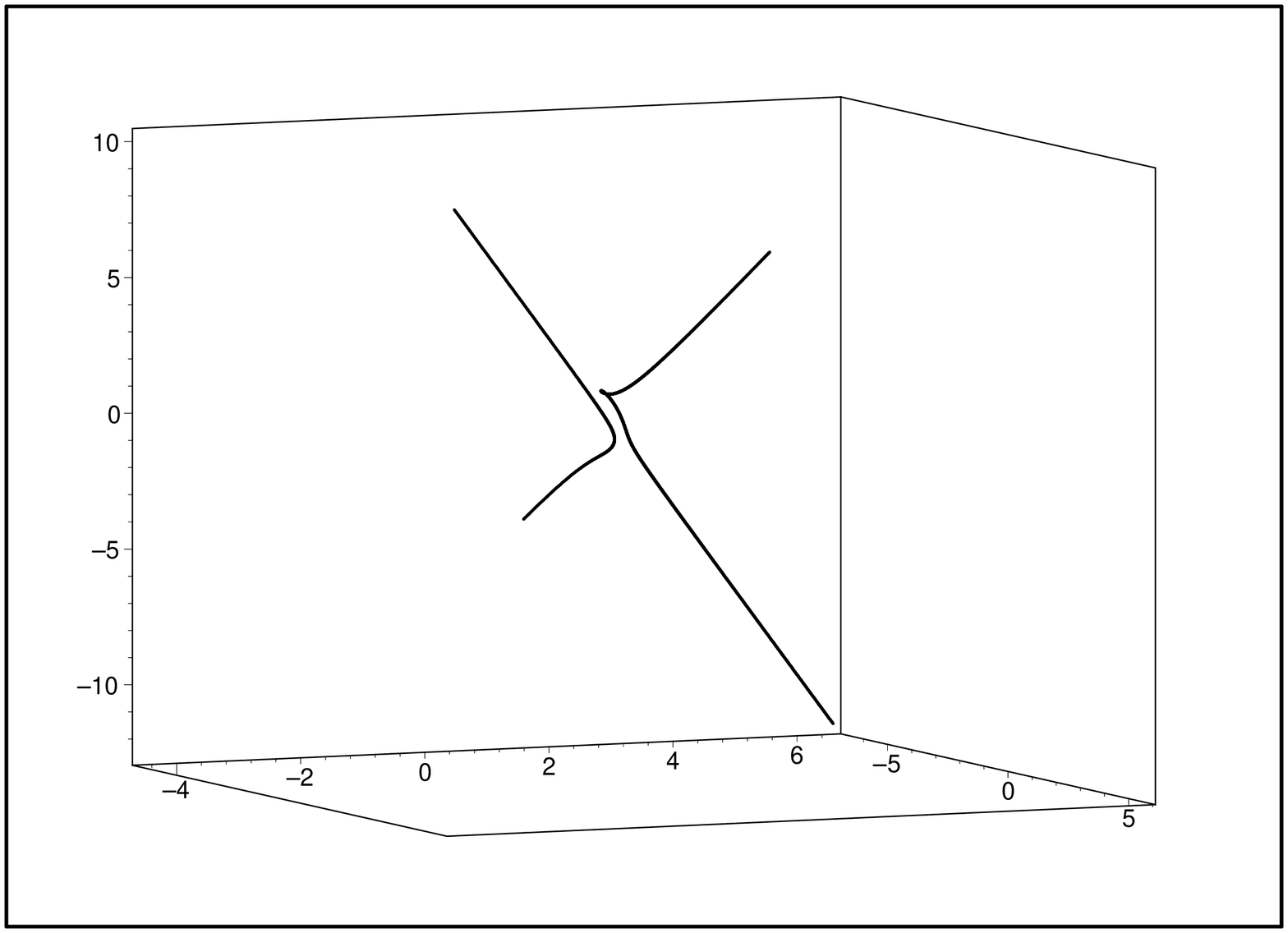,width=6cm,height=6cm,angle=270}} \\
\mbox{\psfig{figure=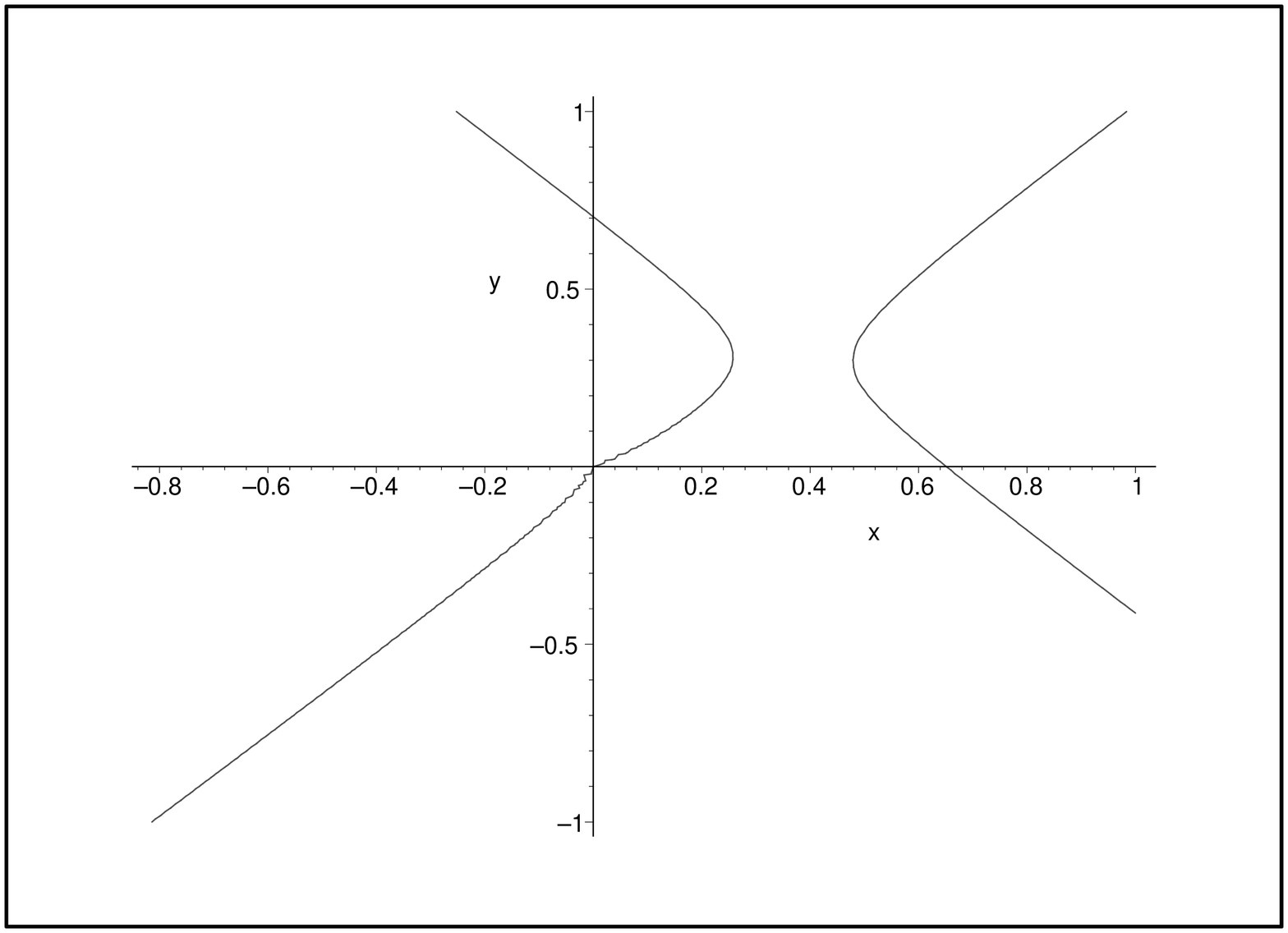,width=6cm,height=6cm,angle=270}} &
\mbox{\psfig{figure=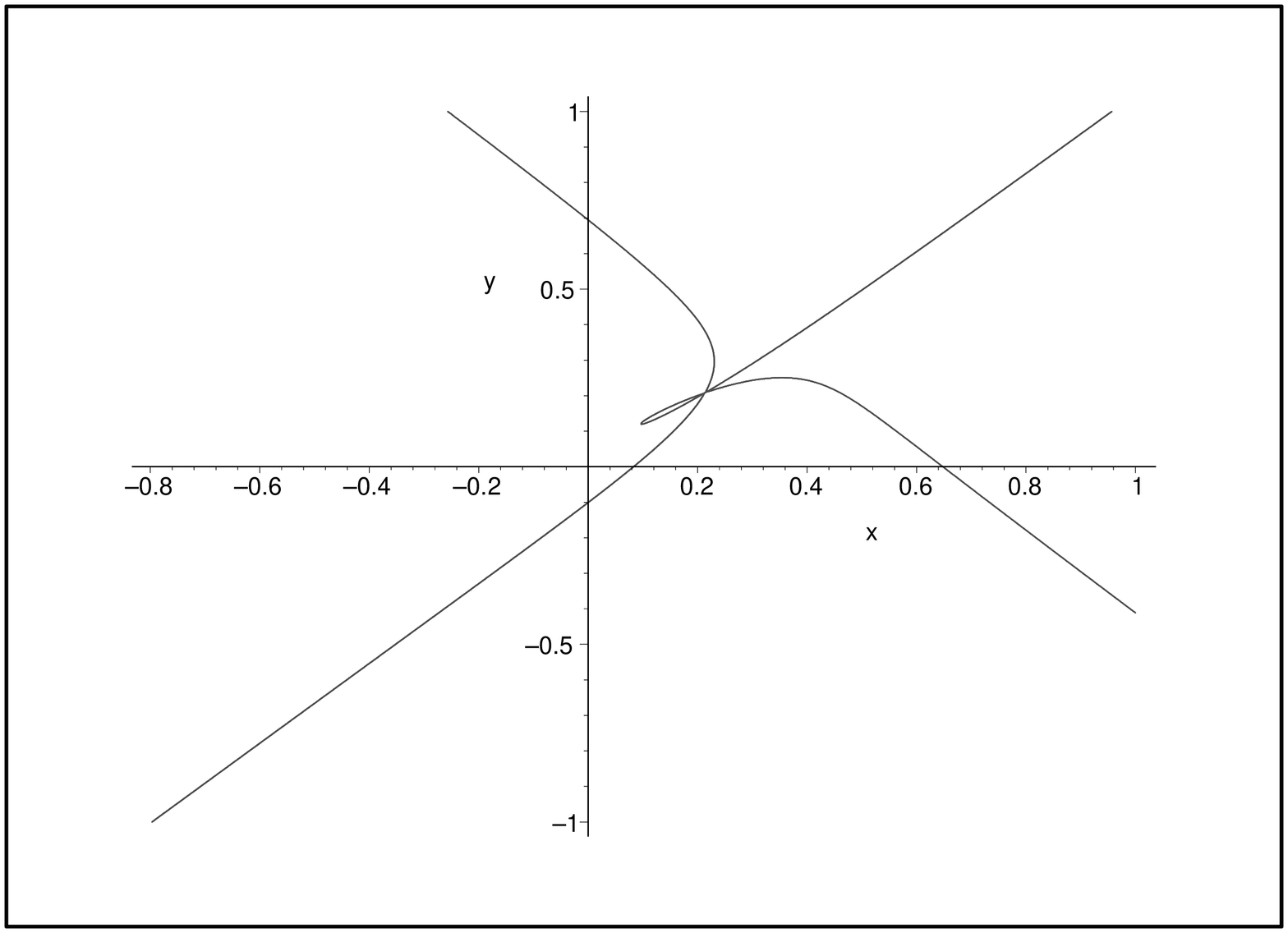,width=6cm,height=6cm,angle=270}} \\
\end{array}$}
\caption{Up left: plot of $\cc$; Up right: plot of $\overline{\cc}$; Down left: plot of $\pi_z(\cc)$; Down right: plot $\cal D$.}
\end{figure}\label{fig-1}
\end{center}

It only remains to compute the numerator of the third component of the rational parametrization of the lifted curve; namely $p_3(t)$. Applying Lemma \ref{lemma-interpolation} we get the approximate   parametrization (see Fig. \ref{fig-conjunta}  $\cP (t)=(\frac{p_1(t)}{q(t)},\frac{p_2(t)}{q(t)}, \frac{p_3(t)}{q(t)})$ of the space curve $\cc$ (i.e. the parametrization of $\overline{\cc}$),  where
\[ \frac{p_{3}(t)}{q(t)}={\frac {- 1.067157288\,{t}^{3}- 0.2783759249- 0.7182447737\,t+
 1.955832944\,{t}^{2}}{- 0.9059858774+ 1.956830479\,t- 0.6103552658\,{
t}^{2}- 1.494650450\,{t}^{3}+{t}^{4}}}.
 \]
\begin{center}
\begin{figure}[ht]
\centerline{
\psfig{figure=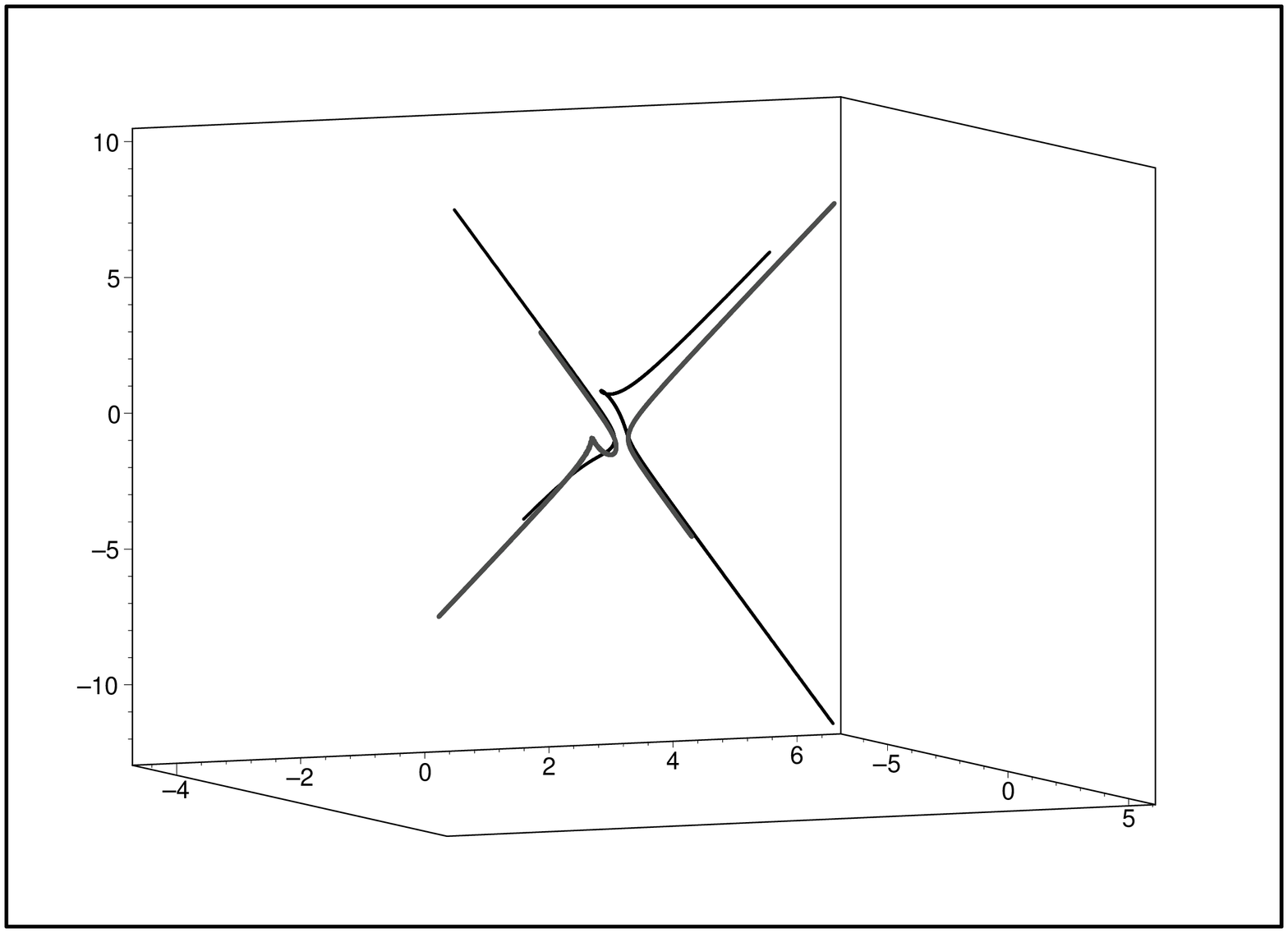,width=10cm,height=10cm,angle=270}}
\caption{Joint plot of $\cc$ and $\overline{\cc}$}
\end{figure}\label{fig-conjunta}
\end{center}

\end{example}

\begin{example}\label{ex-2}
Let $\cc$ be the space curve defined by the polynomials
\begin{align*}
\mathit{F_1} =& 20052827033xy+2850904342zy-7155364672zx-\frac{215763180597}{100}x\\ & -7869010116z+ \frac{1743412651801}{100}y-43102722226y^2+1610946062z^2 \\
\mathit{F_2} =& -18330943984zy+33857630124zx-\frac{390188402999}{25}x-56921602320z+ \\ & \frac{12611223036001}{100}y-166608514760y^2+57179742076z^2+20052827033x^2
\end{align*}
and let $\epsilon=\frac{1}{600}$. One  can check that $\cc$ satisfies all the hypotheses imposed in Section \ref{sec-general-assumptions}. Moreover, $\deg(\cc)=4$. The projected curve  $\pi_z (\cc)^*$ is not $\epsilon$-rational, but $\pi_y(\cc)$ is. So we work with the projection on the plane $y=0$. The defining polynomial of $\pi_y (\cc)^*$ is
\para

\noindent
\begin{flushleft}
$
f(x,z)={ 6.959832072\cdot 10^{47}}\,zx-{ 4.075715387\cdot 10^{36}}\,x-{
 6.207866771\cdot 10^{35}}\,z+{ 1.769623619\cdot 10^{47}}\,{x}^{3}+{
 9.541705261\cdot 10^{46}}\,{x}^{2}+{ 1.269145848\cdot 10^{48}}\,{z}
^{2}+{ 8.077561390\cdot 10^{47}}\,{x}^{3}z-{ 1.355904241\cdot 10^{48
}}\,{x}^{2}z-{ 2.573514563\cdot 10^{48}}\,{z}^{3}+{ 2.289865008
\cdot 10^{48}}\,{z}^{4}-{ 3.292700550\cdot 10^{48}}\,{z}^{2}x+{
 4.798217962\cdot 10^{48}}\,x{z}^{3}+{ 3.090311649\cdot 10^{48}}\,{x
}^{2}{z}^{2}-{ 2.981944666\cdot 10^{47}}\,{x}^{4}.
$
\end{flushleft}

\para

Note that $\deg(\pi_y (\cc)^*)=4$. Applying the
 approximate parametrization algorithm for plane curves in \cite{PRSS} we get that  $\pi_y (\cc)^*$ is $\epsilon$-rational.
  Furthermore  the algorithm outputs $\cQ (t)=(\frac{p_1(t)}{q(t)},\frac{p_3(t)}{q(t)})$
\begin{align*}
p_1=& -0.8393730698 t-1.927074242 t^2-0.1348623847-0.7172362732t^4-1.9338152490t^3, \\
p_3=& 0.009163725887 t-0.08721680470
t^3+0.004664930306-0.05315337574 t^4
\\
& -0.03060855191 t^2,
 \\
q= &  3.555348439\,{t}^{3}+{t}^{4}+ 4.622830832\,{t}^{2}+
2.625458073\,t+
 0.5529230644.
\end{align*}
It only remains to compute the numerator of the second component of
the rational parametrization of the lifted curve; namely $p_2(t)$.
Applying {\sf Algorithm-Lift} we get the    parametrization  $\cP
(t)=(\frac{p_1(t)}{q(t)},\frac{p_2(t)}{q(t)}, \frac{p_3(t)}{q(t)})$
of the space curve $\ccc$  where
\begin{align*}
p_2=& 0.3088642466 t+0.4323589964 t^2+0.07398086924+0.2023592378
t^3.
\end{align*}

\end{example}


\section{Error Analysis}\label{sec-error-analysis}

In this section, we prove that the Hausdorff distance between the input and output curves, of our algorithm, is always finite.  For this purpose, we first need to develop some results on asymptotes of space curves. Afterwards we will analyze the distance. To start,  we briefly recall the notion of Hausdorff distance; for further details we refer to \cite{AB}. In a metric space $(X,\ddd)$,  for $\emptyset \neq B\subset X$ and $a\in X$ we define
\[ \ddd(a,B)=\inf_{b\in B}\{\ddd(a,b)\}. \]
Moreover,  for  $A,B\subset X\setminus \emptyset$ we define
\[ \HH_\ddd(A,B)=\max\{\sup_{a\in A}\{\ddd(a,B)\}, \sup_{b\in B}\{\ddd(b,A)\} \}.  \]
By convention $\HH_\ddd(\emptyset,\emptyset)=0$ and, for $\emptyset \neq A\subset X$, $\HH_\ddd(A,\emptyset)=\infty$. The function $\HH_\ddd$ is called the {\sf Hausdorff distance induced by $\ddd$}. In  our case, since we will be working in $(\C^3,\ddd)$ or $(\R^3,\ddd)$, being $\ddd$ the usual unitary or Euclidean distance, we simplify the notation writing $\HH(A,B)$.

 \para

Let $\ee$ be an space curve in $\C^3$; similarly if we consider the curve in $\C^n$. The intuitive idea of asymptote is clear, but here we formalize it and state some results.
 Although these results might be  part of  the background on  the  theory of asymptotes, we have not been able to find a  reference in the literature suitable for our needs.

  \para

  We say that a line $\cal L$ in $\C^3$ is an {\sf asymptote} of  $\ee$ if there exists a sequence $\{P_n\}_{n\in \mathbb{N}}$ of points in $\ee$ such that
$\lim_{n} \| P_n \|=\infty$ and $\lim_n \ddd(P_n, {\cal L})=0$; where $d$ denotes the usual unitary distance in $\C^3$ and $\|\,\|$ the associated norm.  In  the following, we show how the tangents at the simple points at infinity of $\ee$ are related with the asymptotes.
More precisely, we have the following lemma. In the sequel, if $\cal V$ is a projective variety in $\proj^3(\F)$, we denote by ${\cal V}_a$ the open set  ${\cal V} \cap \{w\neq 0\}$. In addition, for $\lambda \in  \C$, $\overline{\lambda}$ denotes its conjugate.

\para

\begin{lemma}\label{lemma-asintotas}
Let $P=(a:b:c:0)\in \ee^\infty$ be simple and let $\tgh$  be the tangent line to $\ee^h$ at $P$. If $\tgh$ is not included in the plane $w=0$, then $\tgh_a$ is an asymptote of $\ee$. Moreover, $(\overline{a},\overline{b},\overline{c})$ is a direction vector of the asymptote.
\end{lemma}

\para

\noindent {\bf Proof.} Let $\{H_1,\ldots,H_m\}\subset \F[x,y,z,w]$ be homogeneous polynomials defining the ideal of  $\ee^h$. Let
\[ \pi_{i}(x,y,z,w)= \frac{\partial H_{i}}{\partial x}(P) x+ \frac{\partial H_{i}}{\partial y}(P) y+
\frac{\partial H_{i}}{\partial z}(P) z+ \frac{\partial H_{i}}{\partial w}(P) w. \]
Then, $\tgh$ is the projective variety defined by $\{\pi_1(x,y,z,w),\ldots,\pi_m(x,y,z,w)\}$
(see e.g. pp. 181 in  \cite{harris}).
Note that, since $\tgh$ is  not included in the plane $w=0$,  $\tgh_a$ is an affine line. Now, we consider a local parametrization $\cP(t)=(\tilde{x}(t):\tilde{y}(t):\tilde{z}(t):w(t))$ of $\ee^h$ centered at $P$. Since  $P$ is simple, and its tangent is not included in $w=0$, the multiplicity of intersection of $\ee^h$ and the plane $w=0$ at $P$ is 1. Thus, $w(t)$ can be expressed as  $w(t)=tu(t)$, where $u(t)$ has order 0. Therefore, $u(t)^{-1}$ is a  power series and $\cP(t)$ can be expressed as
\[ \cP(t)=(x(t):y(t):z(t):t) \]
where $x(t)=\tilde{x}(t) u^{-1}(t)$, $y(t)=\tilde{y}(t) u^{-1}(t)$, $z(t)=\tilde{z}(t) u^{-1}(t)$ are power series. Moreover, since $\tgh$ is the tangent, the order of $\pi_i(\cP(t))$ has to be, at least, 2. Therefore,
 $\pi_i(\cP(t))$ can be expressed as $\pi_i(\cP(t))=t^{\ell_i} v(t),$
where $\ell_i>1$ and $v(t)$ has order 0.

 Now, let  $ \{t_n\}$ be a sequence of complex numbers converging to $0$, and such that $t_n\neq 0$. Then, for all $n$, $P_n=(\frac{x(t_n)}{t_n}, \frac{y(t_n)}{t_n},\frac{z(t_n)}{t_n})\in \ee$. Moreover, $\lim_n \| P_n\|=\infty$. Let $\Pi_i$ be the affine plane defined by $\pi_i(x,y,z,1)$. We prove that, for all $i$,  $\lim_{n} \ddd(P_n,\Pi_{i})=0$. From where, one deduces that $\lim_n \ddd(P_n,\tgh_a)=0$, and hence that $\tgh_a$ is an asymptote.

 \para

 \noindent Indeed,
\[ \ddd(P_n,\Pi_{i})=\frac{\left|\pi_{i}(P_n,1)\right|}{\left\|\left( \frac{\partial H_{i}}{\partial x}(P), \frac{\partial H_{i}}{\partial y}(P), \frac{\partial H_{i}}{\partial z}(P) \right) \right\|}=
\frac{\left| t_{n}^{\ell_i-1} v(t_n)\right|}{\left\|\left( \frac{\partial H_{i}}{\partial x}(P), \frac{\partial H_{i}}{\partial y}(P), \frac{\partial H_{i}}{\partial z}(P) \right) \right\|}.
\]
Since the denominator is not zero, because $\tgh$ is not included in $w=0$, and since $\ell_i-1>0$ and order of $v$ is 0, we get that
$\lim_n \ddd(P_n,\Pi_{i})=0$.

Finally, we see that $u=(\overline{a},\overline{b},\overline{c})$ is a direction vector of the asymptote. First we observe that $u$  is not zero, since $P\in {\Bbb P}^{3}(\F)$. Now, by Euler equality we have that
\[ \frac{\partial H_i}{\partial x} x+ \frac{\partial H_i}{\partial y} y+ \frac{\partial H_i}{\partial z} z+ \frac{\partial H_i}{\partial w} w=\deg(H_i) H_i.\]
Substituting by $P$, we get
\[ \frac{\partial H_i}{\partial x}(P) a+ \frac{\partial H_i}{\partial y}(P) b+ \frac{\partial H_i}{\partial z}(P) c+ \frac{\partial H_i}{\partial w}(P) 0=\deg(H_i) H_i(P)=0.\]
Hence,
\[
\left( \frac{\partial H_i}{\partial x}(P),\frac{\partial H_i}{\partial y}(P),\frac{\partial H_i}{\partial z}(P)
\right) \cdot u =  \frac{\partial H_i}{\partial x}(P) a+ \frac{\partial H_i}{\partial y}(P) b+ \frac{\partial H_i}{\partial z}(P) c=0.
\]
Thus $u$ is orthogonal to all the planes $\Pi_i$. \qed
\para

\begin{remark}\label{remak-asintotas-reales} Given $P$ verifying  the conditions of Lemma \ref{lemma-asintotas}, we observe that:
\begin{enumerate}
\item   The asymptote determined by $P$ can only be real if $P$ is real. This is due to the fact that $\ee^h$ is defined by  real polynomials
\item If the real asymptote $\cal L$  determined by $P$ is real, then there exists a real branch of $\ee$ approaching $\cal L$ in each of the two directions of $\cal L$.
Indeed let $ \cP(t)=(x(t):y(t):z(t):t)$ be, as in the proof of Lemma
\ref{lemma-asintotas}, a real local parametrization centered at $P$.
$x(t), y(t), z(t)$ are power series and at least one of them (say
w.l.o.g. $x(t)$) is of order $0$. Thus the  affine branch tending to
$P$ is locally parametrized by
$\cP^a(t)=\frac{1}{t}(x(t),y(t),z(t))$. Therefore,
\[ \frac{x(t)}{t}=\frac{a_0+a_1t+a_2 t^2+\cdots}{t}=\frac{a_0}{t}+a_1+ t(a_2+\cdots)\,\,\, \mbox{with}\,\, a=a_0\neq 0. \]
Let us assume w.l.o.g. that $a_0>0$. On the other hand
$\lim_{t\rightarrow 0} \frac{x(t)}{t}= \infty$. Therefore, for
positive values $t^+$ of $t$, sufficiently close to $0$, the first
coordinate of $\cP^{a}(t^+)$ is positive, while for negative values
$t^-$ of $t$, sufficiently close to $0$, the first coordinate of
$\cP^{a}(t^-)$ is negative. Thus, the real branch is approaching
infinity following both directions of the asymptote.
\item Let all points in $\ee^\infty$ satisfy  the conditions of Lemma \ref{lemma-asintotas}. Then, $\ee \cap \R^3$ is unbounded iff there is at least one real point in $\ee^\infty$.
\end{enumerate}
\end{remark}

\para



\para

\begin{lemma}\label{lema-asintotas-2}
Let $\ee_1,\ee_2\subset \C^3$  be irreducible real curves  such that
 $\ee_{1}^{\infty}=\ee_{2}^{\infty}$ and  $\card(\ee_{1}^{\infty})=\card(\ee_{2}^{\infty})=\deg(\ee_1)=\deg(\ee_2)$.
Then,
\begin{enumerate}
\item  $\ee_i$ has $\deg(\ee_i)$ different asymptotes, non of them parallel.
\item For each  asymptote $\cal L$ of $\ee_1$ there exists exactly one asymptote $\overline{\cal L}$ of $\ee_2$ such that ${\cal L}$ and $\overline{\cal L}$ are parallel.
\item For each real asymptote $\cal L$ of $\ee_1$ there exists exactly one real asymptote $\overline{\cal L}$ of $\ee_2$ such that ${\cal L}$ and $\overline{\cal L}$ are parallel.
\end{enumerate}
\end{lemma}

\para

\noindent {\bf Proof.} Let $d=\deg(\ee_1)=\deg(\ee_2)$. Since  all
points at infinity of both curves are simple and no tangent line is
included in the plane at infinity,  Lemma \ref{lemma-asintotas}
implies that each curve has exactly $d$ asymptotes. Furthermore,
since $\card(\ee_{i}^{\infty})=d$, again by Lemma
\ref{lemma-asintotas}, the direction vectors are different and not
parallel. This proves (1).  (2) follows from (1) and
$\ee_{1}^{\infty}=\ee_{2}^{\infty}$.  (3) follows from  (2) and
Remark \ref{remak-asintotas-reales} (1). \qed

\para

Applying the previous lemmas we get the next theorem.

\para
\begin{theorem}\label{theorem-distance}
 Let $\ee_1,\ee_2\subset \C^3$  be irreducible real curves  such that
 $\ee_{1}^{\infty}=\ee_{2}^{\infty}$
and
$\card(\ee_{1}^{\infty})=\card(\ee_{2}^{\infty})=\deg(\ee_1)=\deg(\ee_2)$.
Then, $\HH(\ee_1\cap \R^3,\ee_2\cap \R^3)< \infty$.
\end{theorem}

\para

\noindent {\bf Proof.} By Remark \ref{remak-asintotas-reales} (3),
$\ee_1\cap \R^3$ is bounded iff $\ee_2\cap \R^3$ is bounded.
Moreover, if both are bounded the result follows from Lemma 3.58 in
\cite{AB}. If they are not bounded, let
  ${\cal L}_1,\ldots,{\cal L}_n$ be the real
asymptotes of $\ee_1$ and let $\overline{\cal
L}_1,\ldots,\overline{\cal L}_n$ be the  real asymptotes of $\ee_2$,
with ${\cal L}_i\parallel \overline{\cal L}_i$; see Lemma
\ref{lema-asintotas-2}. In addition let $\delta_i$ be the real
branch of $\ee_1$ approaching the asymptote ${\cal L}_i$; similarly
for $\overline{\delta}_i$ and $\overline{\cal L}_i$.  Let
$\gamma=\max\{\HH({\cal L}_i,\overline{\cal
L}_i)\,|\,i=1,\ldots,n\}$; by Lemma 3.58, in \cite{AB}, $\gamma\in
\R$. Now, let $\rho>0$, and consider the cylinder $V_i$ around
${\cal L}_i$ at distance $\rho$; similarly  $\overline{V}_i$ is the
cylinder around  $\overline{\cal L}_i$ at distance $\rho$.

 Let us see that $\ee_1 \nsubseteq (V_1\cup \cdots \cup V_n)$. If $\ee_1$ was included, since $\ee_1$ is irreducible then it would be included in one  $V_i$; say in $V_1$. Then, for every sequence  $\{P_n\}_{n\in \N}$ of points of $\ee_1$ such that $\lim_n\|P_n\|_2=\infty$, one has that $\lim_n \ddd(P_n,{\cal L}_1)=\rho>0$ and $\lim_n \ddd(P_n,{\cal L}_j)=\infty$ for $j\neq 1$ (note that ${\cal L}_1 \not\,\parallel {\cal L}_j$ if $j>1$). Thus, $\ee_1$ would have no  asymptote. Similarly, one has that $\ee_2 \nsubseteq (\overline{V}_1\cup \cdots \cup \overline{V}_n)$.

 Let $\Delta_1=\ee_1 \cap (V_1\cup \cdots \cup V_n)$ and
$\Delta_2=\ee_2 \cap (\overline{V}_1\cup \cdots \cup
\overline{V}_n)$. Because of the reasoning above
$\card(\Delta_1\cup\Delta_2)<\infty$.
 So, let $\alpha=\max\{\|P\|_2 \,|\, P\in \Delta_{1}\cup \Delta_2\}$. In this situation, we can choose a compact ball $B$ centered at the origin, and with a  sufficiently big radius (and bigger than $\alpha$) such that $(\ee_1\cup \ee_2)\cap \R^3$ is included in $B$ union the inner part of the cylinders
  $V_1,\ldots,V_n,\overline{V}_1,\ldots,\overline{V}_n$. Now, let $P\in \ee_1\cap \R^3$, then
  \begin{itemize}
  \item If $P\in B$ then $\ddd(P,\ee_2\cap \R^3)\leq \ddd(P,\ee_2\cap \R^3\cap B)\leq \HH(\ee_1\cap \R^3\cap B,\ee_2 \cap \R^3\cap B)$ that is equal to a real number $M$, because of Lemma 3.58 in \cite{AB}.
 \item If $P\not\in B$ then $P$ is included  in the interior of one of the cylinders $V_i$. Say $P$ is in the interior of $V_1$; so $P\in \delta_1$. By Remark \ref{remak-asintotas-reales} (2), we know that in $\delta_1$ we can distinguish two branches  each in each direction of the asymptote, let us call them $\delta_{1}^{+}$ and $\delta_{1}^{-}$, and say that $P\in \delta_{1}^{+}$; similarly for $\overline{\delta}_{1}, \overline{\delta}_{1}^{+}, \overline{\delta}_{1}^{-}$. In this situation, there exist $A\in {\cal L}_1, B\in \overline{\cal L}_1$ and $Q\in \overline{\delta}_{1}^{+}\cap (\R^3\setminus B)$ such that
     \[ \ddd(P,A)\leq \rho, \ddd(A,B)\leq \gamma, \ddd(B,Q)\leq \rho. \]
     Therefore, $\ddd(P,Q)\leq 2\rho+\gamma$. That is, for every $P\in \ee_1\cap (\R^3\setminus B)$ there exists $Q\in \ee_2\cap \R^3$ such that $\ddd(P,Q)<2\rho+\gamma$; so $\ddd(P,\ee_2\cap \R^3)<2\rho+\gamma$.
  \end{itemize}
  Since the reasoning above can be done analogously for $Q\in \ee_2\cap \R^3$, we conclude that
  \[ \HH(\ee_1\cap \R^3,\ee_2\cap \R^3)\leq \max\{M,2\rho+\gamma\}. \]

\para

\begin{corollary}\label{cor-distance}
 Let $\cc, \ccc$ be the input and output curves of our algorithm. Then $\HH(\cc\cap \R^3,\ccc\cap \R^3)< \infty$.
\end{corollary}

\para

\noindent {\bf Proof.} It follows from the general assumptions, and from Theorems \ref{teorema-lifted-curve} and \ref{theorem-distance}. \qed

\end{document}